\documentclass[10pt]{amsart}
\DeclareSymbolFont{AMSb}{U}{msb}{m}{n}

\usepackage{amssymb}
\usepackage{amsfonts}
\usepackage{amsmath}
\usepackage{amsthm}

\usepackage[english]{babel}
\usepackage{graphicx}

\newtheorem*{thm}{Theorem}
\newtheorem{lem}{Lemma}
\newtheorem*{prop}{Proposition}

\begin{document}

\title[D. Gage \& D. Jackson]{Computer Generated Images for Quadratic Rational Maps with a Periodic Critical Point}

\author{Dustin Gage}
\address{Division of Mathematics and Computer Science, University of Maine at Farmington, Farmington, ME 04938}
\email{dustin.gage@maine.edu}
\urladdr{http://students.umf.maine.edu/dustin.gage/public.www/}

\author{Daniel Jackson}
\address{Division of Mathematics and Computer Science, University of Maine at Farmington, Farmington, ME 04938}
\email{daniel.jackson1@maine.edu}
\urladdr{http://faculty.umf.maine.edu/daniel.jackson1/public.www/}

\thanks{During a portion of the work on this paper, the first author was supported by two University of Maine at Farmington, Michael D. Wilson Scholarships for undergraduate research.}

\keywords{rational map, complex dynamics, plane curve singularities, geometric genus, hyperbolic maps, Mandelbrot set}
\subjclass{37F45,37F10,14H50}
\date{September 17, 2010}

\pagestyle{plain}

\begin{abstract} We describe an algorithm for distinguishing hyperbolic components in the parameter space of quadratic rational maps with a periodic critical point. We then illustrate computer images of the hyperbolic components of the parameter spaces $V_1-V_4$, which were produced using our algorithm. We also resolve the singularities of the projective closure of $V_5$ by blowups, giving an alternative proof that as an algebraic curve, the geometric genus of $V_5$ is 1. This explains why we are unable to produce an image for $V_5$. 
\end{abstract}

\maketitle

\section{Introduction}

$V_n$ is the set of holomorphic conjugacy classes of quadratic rational maps with a critical point of period $n$. For example, $V_1$ may be identified with the family of quadratic polynomials $$z \mapsto z^2+c, \textrm{ }c\in \mathbb C,$$ each map having $\infty$ as a fixed critical point. The dynamics of the maps in $V_1$ are encoded by the much studied Mandelbrot set (see e.g. [BM],[M],[DH]). 

$V_2$ can be taken as the family of functions $f_a(z)=\frac{a}{z^2+2z}$ (along with the function $f(z)=1/z^2$), each map having the critical 2-cycle $0\mapsto \infty \mapsto 0.$ In [T], Timorin has given a detailed description of the dynamical plane of $V_2$. 

In the case of $n\geq3$, one can take $V_n$ to be the set of $(b,c)\in \mathbb C^2$ such that for the quadratic map $$f(z)=1+b/z+c/z^2,$$ the critical point $0$ has period $n$ (see e.g. [R3]). So, for $n\geq3$, any map $f$ in $V_n$ has the critical cycle $$0 \mapsto \infty \mapsto 1 \mapsto ... \mapsto f^{n-2}(1)=0.$$ 

$f^{n-2}(1)=0$ may be written in the form $$\frac{P_n}{Q_n}=0,$$ where $P_n,Q_n \in \mathbb C[b,c]$ are polynomials having no common factors. The closure of $V_n$ (in $\mathbb C^2$) is a complex algebraic curve supported on $P_n=0$. 

For example, since $f(1)=1+b+c$, $V_3$ identifies with the complex line $$P_3=1+b+c=0.$$ 

For $n=4$, we have $$f^{2}(1)=\frac{1+3b+2b^2+3c+3bc+c^2}{P_3^2}.$$ So $V_4$ is contained in the zero set of the irreducible conic $$P_4=1+3b+2b^2+3c+3bc+c^2.$$ 

Now, for $1 \leq n \leq 4$, the projective closure of $V_n$, which we will denote by $\overline{V_n}$, is birational to $\mathbb P^1=\mathbb C \cup \{\infty\}.$ However, as is proved in Stimson's thesis [St], and as we will demonstrate in section 5, the geometric genus of $\overline{V_5}$ is 1, i.e. $\overline{V_5}$ is birational to a torus. In [St] it is also shown that $g(V_6)=6$ and conjectured that $g(V_7)=22$ (with the irreducibility of $V_7$ left unfinished). This means that in these higher period cases, we no longer have a parametrization of $V_n$ by $\mathbb C$.

A rational map $f:\mathbb P^1=\mathbb C \cup \{\infty\} \rightarrow \mathbb P^1$ is hyperbolic if under iteration, each critical point of $f$ is attracted to some attracting periodic cycle. Hyperbolic maps are an open and (conjecturally) dense set in the space of rational maps (see e.g. [MSS]). Connected components of the set of hyperbolic maps are called hyperbolic components. Over the past several decades, much progress has been made in describing the hyperbolic components of the $V_n$'s - especially for $n\leq3$ (see e.g. [R1-3],[W],[T]). 

In this paper we will describe an algorithm to draw the hyperbolic components of $V_n$ while distinguishing between types of components. Our algorithm is essentially an implementation of the classification of hyperbolic components for rational maps (see e.g. [R1-3]). We then use this algorithm to generate graphical approximations of the components of $V_n$ for the cases for which we have a parametrization (the genus 0 cases of $V_1-V_4$). Of course $V_1-V_3$ have been drawn before (see e.g. [BM],[T],[W]), but we have included images of these parameter spaces for the reader's convenience (see Figs. 2,3, and 4). Our graphical approximations of $V_4$ are shown in Figs. 5,6, and 7. Recently, Kiwi and Rees ([KR]) have found fomulas to count the number of hyperbolic components in $V_n$. Their results confirm the number of components illustrated in our figures. 

To complete the paper, we give an additional proof that the genus of $V_5$ is 1 (Theorem 1). We do our genus calculation for $V_5$ by resolving singularities by blowups, whereas the original proof in [St] uses Puiseaux expansions. This additional proof has been included for 2 reasons. The first reason being completeness of this paper - the result explains why we are not able to generate graphics for $V_5$. The second reason being that the authors hope that a resolution by blowups may shed light on a general genus computation for $V_n$. During the preparation of this paper, the second author and two undergraduate students performed a similar calculation for $V_6$. The blowup sequences for $V_6$ may be found in [DHJ].  

This paper is organized as follows: In section 2 we recall the classification of hyperbolic quadratic rational maps with a critical cycle. We then explain how this classification gives an algorithm for approximating the hyperbolic components of $V_n$, while distinguishing between their types. In section 3 we illustrate some computer generated representations of the hyperbolic components of $V_n$ for $n=1,...,4$. Section 4 contains our calculation of the genus of $V_5$. Finally, section 5 describes the computer program we made to produce the graphics in this paper.

The second author would like to thank Mary Rees for providing some helpful comments and pointing out some useful references during the preparation of this paper.
 
\section{Hyperbolic Quadratic Rational Maps}

Since quadratic maps have 2 critical points and any map $f$ in $V_n$ has the critical cycle $$0 \mapsto \infty \mapsto 1 \mapsto ... \mapsto f^{n-2}(1)=0,$$ we shall refer to the one remaining critical point of $f$ as the free critical point. Hyperbolic rational maps have been classified by their critical orbits (see e.g. [R1-3]). In fact, any hyperbolic map $f$ in $V_n$ must be exactly one of the following four types: 

{\bf Type 1 } The free critical point is attracted to the other critical point, which is a fixed point.

{\bf Type 2 } The free critical point is in a periodic component of the attracting basin of the critical cycle.

{\bf Type 3 } The free critical point is in a preperiodic component of the attracting basin of the critical cycle. 

{\bf Type 4 } The free critical point belongs to the attracting basin of a periodic orbit other than the critical cycle.

This classification suggests an algorithm for making an approximation of the hyperbolic components of $V_n$. Note that it is easy to distinguish type 4 mappings from the other types by testing the orbit of the free critical point for attraction to the critical cycle. Also, if $n>1$ then $V_n$ has no type 1 maps, and if $n=1$ then $V_n$ has no maps of types 2 or 3. So the case of $n=1$ is easy; and, for $n>1$ we just need to distinguish between type 2 and type 3 maps. To do this, one must decide if there is a path from the free critical point to the attracting periodic point, entirely contained within the immediate attracting basin. For most type 2 maps in $V_n$, the line segment between the free critical point and its attractor lies within the immediate basin. However, there are type 2 maps in $V_3$ and $V_4$ for which a nonlinear path must be found (see Figs. 1,7, and 9). We solved this problem by flood-filling the immediate attracting basin to test for the free critical point. 
 
\begin{figure}[!ht]
  \centering    \includegraphics[width=.75\textwidth]{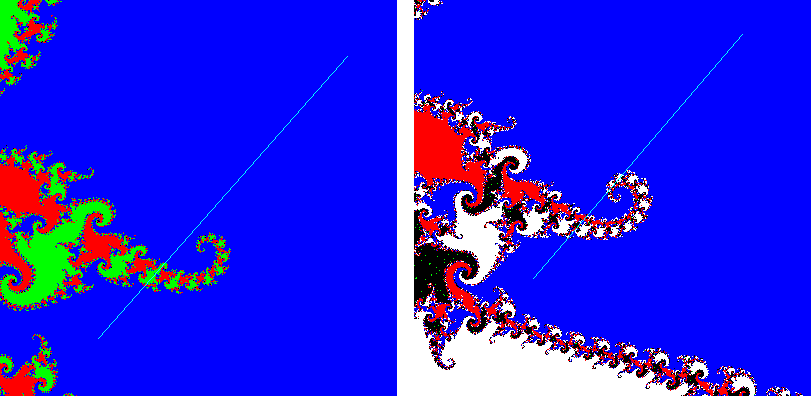}
\caption{The image on the left shows the basins of attraction for the map $f(z)=(z-c)(z-1)/z^2$ where $c=.16-2.2i$. The free critical point is attracted to 1 under iteration of $f^3$ in this case; but, as is illustrated, the line segment joining the free critical point and 1 is not entirely contained in the immediate basin. The image on the right shows a similar situation for a map in $V_4$.}
\end{figure} 

\section{Computer Generated Images}

Using the above algorithm (and a parametrization for $V_n$), we can generate a graphical approximation of the hyperbolic components of $V_n$ - distinguishing between types. In this section we illustrate such approximations for $n=1,2,3,$ and $4$. Approximations of the hyperbolic components of $V_1,V_2,$ and $V_3$ have previously been illustrated (see e.g. [B-M],[M], and [W]). However, these approximations (at least for $V_2$ and $V_3$) do not give a graphical distinction between the types of components. So we will include these parameter spaces for the reader's convenience. Next we will provide some representations of $V_4$.
  
{\bf $V_1$:} The set of holomorphic conjugacy classes of quadratic rational maps with a fixed critical point may be identified with the family of polynomials $$f(z)=z^2+c \textrm{ where } c\in \mathbb C.$$ There are no type II or III components in this case. The complement of the single type I component is the classical Mandelbrot set (see Fig. 2).  

\begin{figure}[!ht]
  
  \centering
    \includegraphics[width=.25\textwidth]{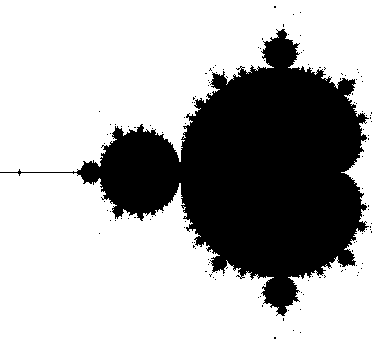}
\caption{$V_1$: The single type 1 component is in white. Since $\infty$ is a fixed critical point there are no maps of types 2 or 3. The free critical point, for the maps colored in black, is not attracted to $\infty$.}

\end{figure}

{\bf $V_2$:} The set of holomorphic conjugacy classes of quadratic rational maps with a period 2 critical point may be identified with the rational map
$$f(z)=\frac{1}{z^2}$$ 
together with the family
$$f(z)=\frac{a}{z^2+2z}, \ a \in \mathbb C - \{0\}.$$

$V_2$ contains one type 2 component (see Fig. 3). A detailed description of the hyperbolic components of $V_2$ has been given in [T].

\begin{figure}[!ht]
\centering
\includegraphics[width=.75\textwidth]{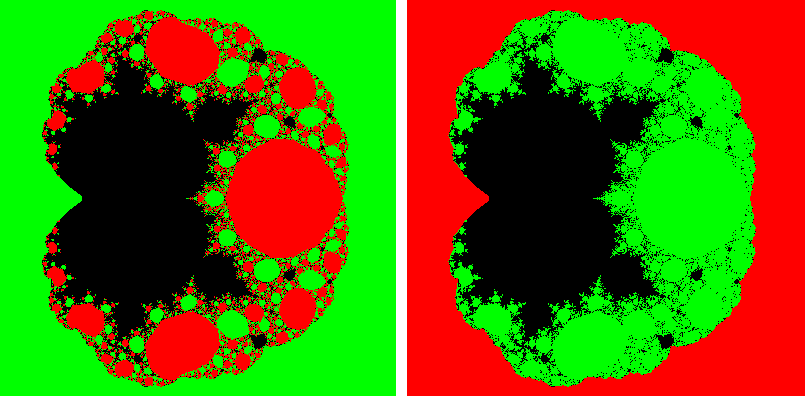}
\caption{$V_2$: The image on the left shows the maps whose free critical point is attracted to $\infty$ in red, while the maps with attraction to 0 are colored green. The image on the right shows the type 2 hyperbolic component in red and the type 3 components in green. In both images the free critical point for the maps colored in black is not attracted to the critical cycle.}
\end{figure}

{\bf $V_n,n \geq 3$:} In this case, $V_n$ is the collection of quadratic functions $$f(z)=1+\frac{b}{z}+\frac{c}{z^2}$$ that have the critical $n$-cycle $$0 \mapsto \infty \mapsto 1 \mapsto ... \mapsto f^{n-2}(1)=0.$$ 

For example, $V_3$ is defined by $$f(1)=1+b+c=0,$$ which gives us the one parameter family of maps $$f(z)=1+\frac{-1-c}{z}+\frac{c}{z^2}=\frac{(z-1)(z-c)}{z^2}.$$ 

$V_3$ has 2 type II components: one containing $c=-1$ and the other containing $c=1$ (see Fig. 4). A nearly complete topological description of the hyperbolic components of $V_3$ has been given in [R3].

\begin{figure}[!ht]
\centering \includegraphics[width=.75\textwidth]{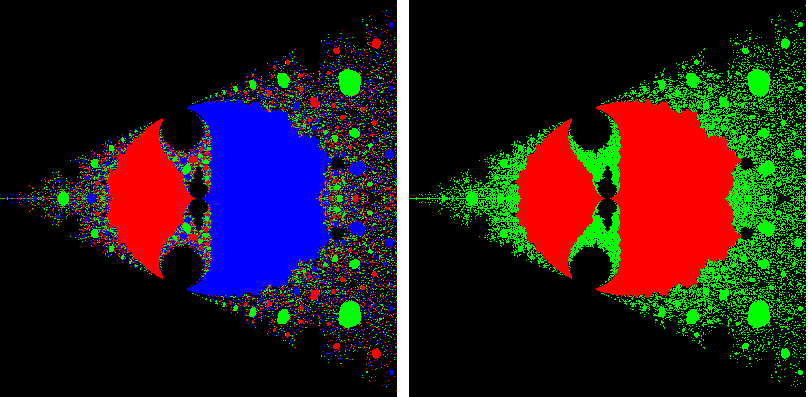} \caption[V3.]{$V_3$: The image on the left shows the maps whose free critical point is attracted to $\infty$, 0, and 1 in red, green, and blue, respectively. The image on the right shows the type 2 maps in red and the type 3 maps in green. In both images the free critical point for the maps colored in black is not attracted to the critical cycle. }\label{fig:fig1}\end{figure}

$V_4$ is determined by $f^{2}(1)=0$, which defines the algebraic curve $$1+3b+2b^2+3c+3bc+c^2=0.$$ One may find a rational parametrization for an irreducible conic by projecting from any point on the curve (see e.g. [S]). Figures 5 and 6 show some representations of $V_4$ using two different parametrizations. Our images show $V_4$ to have 6 type 2 components - which is confirmed by the formulas given in [KR]. 

\begin{figure}[!ht]
	\centering
		\includegraphics[width=.75\textwidth]{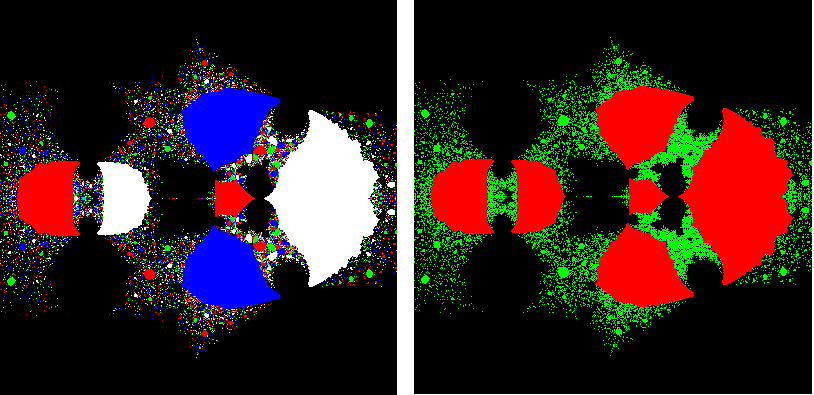}
	\caption[V4.]{$V_4$ projected from $(-1,0)$: The image on the left shows the maps whose free critical point is attracted to $\infty$, 0, 1, and $1+b+c$ in red, green, blue, and white, respectively. The image on the right shows the type 2 components in red, while the green represents the type 3 components. In both images the free critical point for the maps colored in black is not attracted to the critical cycle.}
\end{figure}

\begin{figure}[!ht]
	\centering
		\includegraphics[width=.75\textwidth]{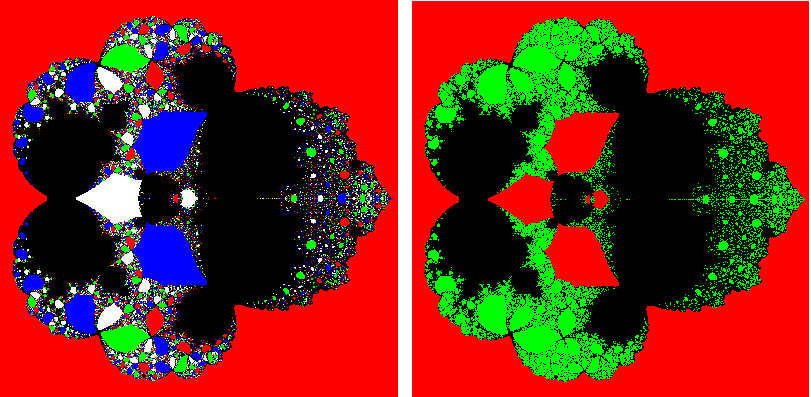}
	\caption[V4.]{$V_4$ projected from $(0,(3+\sqrt{5})/2)$: The coloring is the same as in Fig. 5.}
\end{figure}

Both $V_3$ and $V_4$ required the flood-fill algorithm to distinguish their type 2 components from their type 3 components (see Fig. 7). $V_2$ does not seem to have any such maps.

\begin{figure}[!ht]
	\centering
		\includegraphics[width=.75\textwidth]{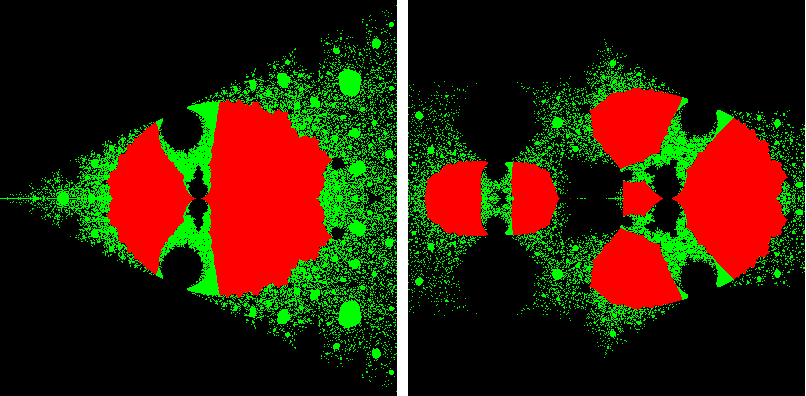}
	\caption[Vn Linear.]{The red approximates those type 2 maps for which the line segment between the free critical point and its attracting periodic point is entirely contained within the immediate attracting basin. $V_3$ is on the left, while $V_4$ is on the right. Compare with Figs. 4 and 5.}
\end{figure}

\section{The Genus of $V_5$}

It is a classical result that any singular point $p$ on an algebraic curve $C$ may be resolved by a sequence of blowups (see e.g. [S] or [C]). The point $p$ and the singular points arising from these blowups are called the infinitely near points to $p$.

To calculate the geometric genus of an irreducible complex projective plane algebraic curve $C\subset \mathbb P^2$ one may use the genus formula:

$$g(C)=\frac{(d-1)(d-2)}{2}-\sum \frac{k_i(k_i-1)}{2},$$ 

where $d$ is the degree of $C$ and the $k_i$'s are the multiplicities of all the infinitely near points to the singularities of $C$ (see e.g. [S]).

Geometrically, blowing-up a point consists of replacing the point by a line of tangent directions. Algebraically, the blowup of a point in affine space $\mathbb A^n = \{(a_1,a_2,...,a_n)|a_i \in \mathbb C\}$ is described as follows:

After suitable change of coordinates one may arrange that the point to blowup is the origin $(0,0,...,0)\in \mathbb A^n$. In this case, setting   

$$B=\{(a_1,a_2,...,a_n;x_1:x_2:...:x_n)|a_ix_j=a_jx_i \textrm{ for }i,j=1,...,n\} \subset \mathbb A^n \times \mathbb P^{n-1},$$

the birational map $\pi:B \rightarrow \mathbb A^n$ given by $$\pi(a_1,a_2,...,a_n;x_1:x_2:...:x_n)=(a_1,a_2,...,a_n),$$ is called the blowup of $\mathbb A^n$ centered at the origin.

For more information on plane curves, singularities, and blowups the reader may refer to a resource such as [S] or [C].

For $n \geq 3$, $V_n$ is contained in a complex plane algebraic curve defined by $$f^{n-2}(1)=0,$$ where  $$f(z)=1+\frac{b}{z}+\frac{c}{z^2}.$$ 

$V_5$ is defined by $$f^3(1)=0,$$ which simplifies to $$\frac{P_5}{P_4^2}=0,$$ where $P_4$ is the polynomial defining $V_4$ (see above) and $P_5 \in \mathbb C[b,c]$ is a degree 5 polynomial. The homogenization of $P_5$ is

\begin{align*}H(a,b,c) = \ &a^5 + 7a^4b + 18a^3b^2 + 21a^2b^3 + 11ab^4 + 2b^5+\\
&7a^4c + 33a^3bc + 53a^2b^2c + 35ab^3c + 8b^4c +\\ 
&15a^3c^2 + 44a^2bc^2 + 42ab^2c^2 + 13b^3c^2 +\\
&12a^2c^3 + 23abc^3 + 11b^2c^3 +\\
&5ac^4 + 5bc^4 + c^5,\end{align*}

(i.e. $P_5=H(1,b,c)$). In what follows we shall refer to the complex projective plane algebraic curve $H=0$ simply as $H$. In the projective coordinates $[a:b:c]$, the singularities of $H$ are at
$$p=[0:1:-1] \textrm{\ \ \ \ and \ \ \ \ } q=[-1:1:0].$$


The infinitely near points to $p$ and $q$ are described by the following 2 lemmas.

\begin{lem} The infinitely near points to $p$ have multiplicities 3, 1, 1, and 1.\end{lem}

\begin{proof} After changing coordinates, we may assume that $p$ occurs at the origin $(a,c)=(0,0)$ on the affine patch $\mathbb A^2_b=\{b=1\} \subset \mathbb P^2$. In this case, the local equation for $H$ is $$a^2c + 3ac^2 + c^3=-(a^5  + 3a^3c + 7a^4c +  8a^2c^2 + 15a^3c^2  + 3ac^3 + 12a^2c^3 + 5ac^4 + c^5).$$ Since $a^2c + 3ac^2 + c^3=c(a+
\frac{3+\sqrt{5}}{2}c)(a+
\frac{3-\sqrt{5}}{2}c)$, $H$ has three distinct tangents at $p$.\end{proof}

\begin{lem}The infinitely near points to $q$ have multiplicities 2, 2, 1, and 1.\end{lem}

\begin{proof}We change coordinates so that $q$ occurs at the origin $(a,c)=(0,0)$ in $\mathbb A^2_b$. Then the local equation for $H$ is $$a^2=g(a,c), \textrm{ where}$$ \begin{align*}
g(a,c)=\ &2a^4 + a^5 - 4a^2c + 5a^3c + 7a^4c - \\&ac^2 - a^2c^2 + 15a^3c^2 - ac^3 + 12a^2c^3 + 5ac^4 + c^5.
\end{align*}

Hence $q$ is a cusp of multiplicity 2.  We will resolve $q$ by blowing up.

Let $\pi:B \rightarrow \mathbb A^2_b$ be the blowup of $\mathbb A^2_b$ at the origin. We can write $B$ explicitly as $$B=\{(a,c;x:y)|ay=cx\} \subset \mathbb A^2_b \times \mathbb P^1,$$ where $[x:y]$ are the projective coordinates on $\mathbb P^1$. Then $H$ is birationally equivalent to the projective closure of

\begin{align*} V =& \ \pi^{-1}( \{a^2= g(a,c)\}-\{(0,0)\}) \cap (\mathbb A^2_b \times \mathbb P^1)\\ =& \ \{(a,c;x:y)|a^2= g(a,c),ay=cx, (a,c) \neq (0,0)\}.\end{align*}

If we set $\infty=(0,0;1:0)$, then using coordinates for the affine patch $\mathbb A^2_b \times \mathbb A^1_y$ we have

\begin{align*} V = \ \{(a,c,x)|&a^2= g(a,c),a=cx, c \neq 0\} \cup \infty\\ = \ \{(a,c,x)|&(cx)^2= g(cx,c),a=cx, c \neq 0\} \cup \infty\\ = \ \{(a,c,x)|& c^2x^2 = c^2(c^3 - cx - c^2x + 5c^3x - 4cx^2 - c^2x^2 + 12c^3x^2 +\\ &5c^2x^3 + 
15c^3x^3 + 2c^2x^4 + 7c^3x^4 + c^3x^5),a=cx, c \neq 0\} \cup \infty \\ = \ \{(a,c,x)|&x(x+c)= c^3  - c^2x + 5c^3x - 4cx^2 - c^2x^2 + 12c^3x^2 +\\ &5c^2x^3 + 15c^3x^3 + 2c^2x^4 + 7c^3x^4 + c^3x^5,a=cx, c \neq 0\} \cup \infty. \end{align*}

This curve in $\mathbb A^3$ has distinct tangent lines $\{a=0,x=0\}$ and $\{a=0,x+c=0\}$ at $(0,0,0)$; therefore, $q$ will be resolved after one more blowup.\end{proof}

\begin{prop}$H$ is irreducible, and hence $\overline{V_5}=H$.\end{prop}

\begin{proof}If $H=C_1 \cup C_2$, then $C_1 \cap C_2$ will be singular points of $H$.  By Bezout's theorem, $C_1 \cap C_2$ consists of $\deg C_1 \times \deg C_2$ points. Counting multiplicities, $H$ has 5 singularities and so $\deg C_1 \times \deg C_2 \leq 5$. Thus $\deg C_1 + \deg C_2 = \deg H = 5$ implies that the degrees of $C_1$ and $C_2$ must be 1 and 4. Let us suppose $C_1$ is a line and $C_2$ is a quartic. Since $C_1 \cap C_2$ consists of 4 points (counting multiplicities) and $H$ has two singularities, $C_1$ must be tangent to $C_2$ at one of the singularities of $H$. From the proof of Lemma 1, the singularity $p$ has 3 distinct tangents and hence $C_1$ is not tangent to $C_2$ at $p$. Thus $C_1$ must be tangent to $C_2$ at the double point $q$. Now, Bezout's theorem implies $C_1$ must intersect $C_2$ two more times. Since $C_1$ is not tangent at $p$, it must meet $C_2$ at a point other than $p$ or $q$. This is contrary to $H$ having only two distinct singular points.\end{proof}

Combining Lemmas 1 and 2 and applying the genus formula given above, we get:

\begin{thm}The geometric genus of $V_5$ is 1.\end{thm}

\section{Generating the Computer Images}

The images in this paper were generated by a Java applet (see Fig. 9) programmed by the authors. Our zoomable fractal generator may be used to graphically explore $V_1,V_2,V_3,$ and $V_4$ (and other spaces of quadratic maps) in a much more detailed manner than given in this paper. The applet is freely available for use and/or download at either authors' websites (see URLs below). The source code for the applet as well as full screen shots of the images in this paper (including input data) are available on the second author's website.

\begin{figure}[!ht]
	\centering
		\includegraphics[width=.75\textwidth]{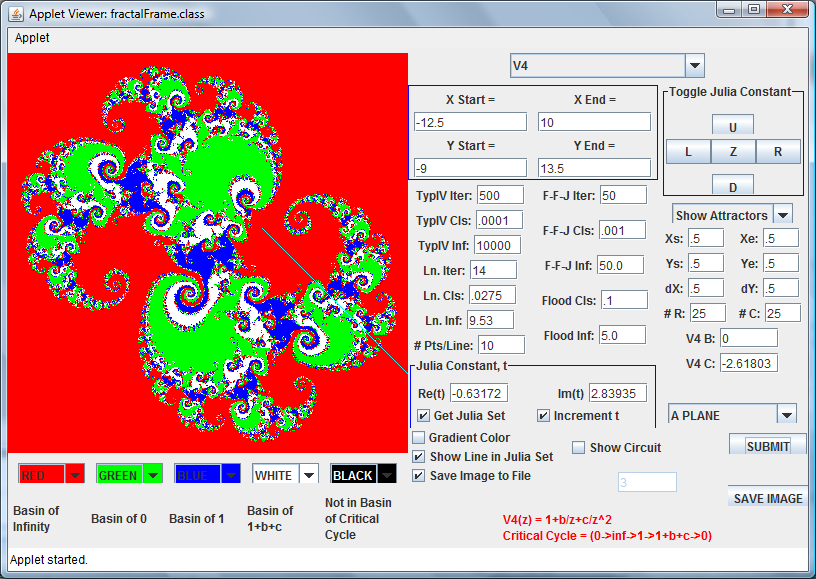}
	\caption[V4.]{Screenshot: The basins of attraction for a type 2 map $f$ in $V_4$. The free critical point is attracted to $\infty$ under iteration of $f^4$. The line segment fails to connect the free critical point with $\infty$.}
	
\end{figure}

\section{References}

[BM] Brooks, R. and Matelski, J. P. {\it Riemann Surfaces and Related Topics: Proceedings of the 1978 Stony Brook Conference. (eds Kra, I. and Maskit, B.)}, Princeton Univ. Press, Princeton, 1981, 65-71.

[C] Casas-Alvero, E. {\it Singularities of Plane Curves}. Cambridge Univ. Press, 2000.

[DHJ] Darby, S., Hall, M., and Jackson, D. {\it A desingularization of $V_6$}. In preparation.


[DH] Douady, A. and Hubbard, J.H. {\it Etude dynamique des polynomes complexes I and II}. Publ. Math. Orsay (1984,1985).

[KR] Kiwi, J. and Rees, M. {\it Counting Hyperbolic Components}. arXiv:1003.6104v1 31 March 2010.

[M] Mandelbrot, B. {\it The Fractal Geometry of Nature.} W.H. Freeman and Company, 1983.

[MSS] Mane, R., Sad,P., and Sullivan, D. {\it On the dynamics of rational maps.} Ann. Sci. Ec. Norm. Sup. {\bf 16} (1983),193-217. 


[Mi] Milnor, J. {\it Geometry and Dynamics of Quadratic Rational Maps}. Experiment. Math. 2 (1993), no. 1, 37--83.

[R1] Rees, M. {\it A Partial Description of the Parameter space of Rational Maps of Degree Two: Part 1}. Acta Math., 168 (1992), 11-87. 

[R2] Rees, M. {\it A Partial Description of the Parameter space of Rational Maps of Degree Two: Part Two}. Proc. Lond. Math. Soc., 70 (1995), 644-690.

[R3] Rees, M. {\it A Fundamental Domain for $V_3$}. Preprint, 2009.


[S] Shafarevich, I. {\it Basic Algebraic Geometry 1.} Second Edition. Springer-Verlag, 1994.

[St] Stimson, J. {\it Degree two rational maps with a periodic critical point.} Thesis, University of Liverpool, 1993. 

[T] Timorin, V. {\it External Boundary of $M_2$}. Fields Institute Communications Vol. 53: Holomorphic Dynamics and Renormalization, A Volume in Honour of John Milnor's 75th Birthday (2008), 225-267.

[W] Wittner, B. {\it On the bifurcation loci of rational maps of degree two}. Ph.D. Thesis, Cornell University, 1988.

\end{document}